\documentclass [a4paper, 10pt]{amsart}
\usepackage[textheight=25 cm, textwidth=16cm, headheight=10pt, headsep=20pt, footskip=10pt]{geometry}

\usepackage{amssymb,mathrsfs}
\usepackage{MnSymbol}
  
\usepackage[hyperindex=true, colorlinks=false]{hyperref}

\newcommand{\cal}{\mathcal }

\newcommand{\R}{{\mathbb R}}

\renewcommand{\a}{\alpha}

\newcommand{\e}{\varepsilon}

\newcommand{\pa}{\partial}

\def\m{{\mu}}

\newtheorem{defi}{Definition}

\def\be{\begin{equation}}
\def\ee{\end{equation}}
\def\bea{\begin{eqnarray}}
\def\eea{\end{eqnarray}}

\def\nn{\nonumber}

\def\o{\omega}

\def\t{\tau}

\def\a {{\alpha}}

\newcommand{\ba}{\begin{aligned}}
\newcommand{\ea}{\end{aligned}}

\newtheorem{Thm}{Theorem}[section]
\newtheorem{Rmk}[Thm]{Remark}
\newtheorem{Prop}[Thm]{Proposition}

\usepackage{color}
\usepackage{hyperref}
\hypersetup{pdfborder=0 0 0}

\newcommand{\norm}[1]{\lVert #1 \rVert}

\newcommand{\bcr}{\begin{color}{red}}
\newcommand{\bcb}{\begin{color}{blue}}
\newcommand{\ec}{\end{color}}
\newcommand{\bcrr}[1]{#1}

\begin{document}

\title {
Mean Field Limit for the Kac model and Grand Canonical Formalism}

\author[T. Paul]{\large Thierry Paul}
\address[T.P] {CNRS  Laboratoire Ypatia des Sciences Math\'ematiques (LYSM),  Roma, Italia \&
Laboratoire Jacques-Louis Lions (LJLL), Sorbonne Universit\'e, Paris, France}  
\email{thierry.paul@sorbonne-universite.fr}

\author[M. Pulvirenti]{\large Mario Pulvirenti}
\address[M.P.] {Sapienza Universit\`a di Roma, Dipartimento di Matematica G.\,Castelnuovo\\ Piazzale A.\,Moro 5, 00185 Roma -- Italia and  International Research Center on the Mathematics and Mechanics of Complex Systems \\ï¿œMeMoCS, University of L'Aquila, Italy}  
\email{pulviren@mat.uniroma1.it}

\author[S. Simonella]{\large Sergio Simonella}
\address[S.S.]{Sapienza Universit\`a di Roma, Dipartimento di Matematica G.\,Castelnuovo\\ Piazzale A.\,Moro 5, 00185 Roma -- Italia}
\email{sergio.simonella@uniroma1.it}

\maketitle

\begin{abstract}

We consider the classical Kac's model for the approximation of the Boltzmann equation, and study the correlation error measuring the defect of propagation of chaos in the mean field limit. This contribution  is inspired by a recent paper of the same authors \cite{PPS} where a large class of models, including quantum systems, are considered. Here we outline the main ideas in the context of grand canonical measures, for which both the evolution equations and the proof simplify.

\end{abstract}

\tableofcontents
\LARGE

{\it This note is dedicated to the memory of our colleague and friend Maria Conceicao Carvalho. Among her many important contributions to the theory of Kac's model, Sao's work \cite{CCG99,CCG05} explained the long time behaviour of the homogeneous Boltzmann equation without using the 
$H$-theorem nor global properties, but exploiting a Dyson expansion over trajectories. 
Our contribution is, in spirit, close to this approach. } 

\bigskip

\section{Introduction}

The kinetic description of particle systems is strongly related to propagation of chaos. This property allows one to
substitute the complex dynamics of a huge number of particles by a single nonlinear partial differential equation
for the 
probability density  of a given particle. More precisely, one adopts a statistical description.  At time zero,
the $N-$particle system is assumed to be ``chaotic" in the sense that each particle is distributed identically and independently from the others,  at least up to an error which is vanishing when $N$ diverges. The dynamics creates correlations and the independence property is lost at any positive time. However,
under suitable scaling limits, the statistical independence of any finite group of particles can be recovered, after averaging, in the limit  $N \to \infty$. As a consequence, a given particle 
evolves according to  an effective kinetic equation. The nature of this dynamics is
determined by the microscopic details of the system and by the regime of physical parameters. Such a mechanism 
works in the formal (in a few cases, rigorous) derivation of the most common kinetic equations.

If the system is described by a symmetric probability measure $W^N(Z_N,t)$ at time $t$, where $Z_N=( z_1, \cdots, z_N)$ is a configuration of the system, being $z_i=(x_i,v_i)$ position and velocity of the particle $i$, the probability density of finding the first $j$ particles in the microscopic state $Z_j$  is given by the marginal
$$
f^N_j(Z_j,t) := \int dz_{j+1} \cdots dz_N W^N(Z_N,t)\;.
$$
In case of a two-body interaction, we have that $f^N_1$ depends on $f^N_2$, which depends on $f^N_3$ and so on, 
 and these marginals are not tensorized. Indeed, the dynamics creates correlations even though at time zero 
it is assumed that
$$
\label{ind}
W^N= (f_0)^{\otimes N}\;.
$$
In many cases, the physical regime of interest can be expressed in terms of suitable scaling limits ($N \to \infty$, together with additional prescriptions), which one can hope to use  to prove that
\be
f^N_j (t) \to f(t)^{\otimes j},
\ee
where $f(t)$ is a solution of the aforementioned kinetic equation.

In three dimensions, fundamental scaling limits for Hamiltonian systems are the following.

\begin{itemize}

\item Mean-field limit for classical systems: $N \to \infty$ and the coupling constant of the two-body interaction is scaled by 
$\frac 1N$. 
Namely, a system of point particles with the mean-field Hamiltonian
$$
H_N=\frac 12 \sum_{i=1}^N v_i^2 +\frac 1{2N} \sum_{i \neq j} \varphi (x_i-x_j)\;,
$$
where $\varphi$ is a two-body potential. 
 The corresponding expected kinetic equation is, in this case, the Vlasov equation.
 \\
\item Low-density or Boltzmann-Grad limit: $N \to \infty$ and the diameter of hard-sphere like particles is $\sim N^{-\frac 12}$.
 The corresponding kinetic equation is the Boltzmann equation. \\
 \item Weak-coupling limit: $N \to \infty$ and $\e \sim N^{-\frac 13}$ is a scale parameter coding the interaction length. The radial interaction potential - say, smooth and decaying - is scaled as
$  \varphi (\frac r {\e})$,
 so that the Hamiltonian of the system is
 $$
H_N=\frac 12 \sum_{i=1}^N v_i^2 + \frac{\sqrt{\e}}{2} \sum_{i \neq j} \varphi \left(\frac {|x_i-x_j|}{\e}\right).
$$
 The corresponding kinetic equation is the Landau equation. \\
 \item High-density limit: $N \to \infty$ and $\e \sim N^{-\frac 14}$, where the Hamiltonian is
 $$
H_N=\frac 12 \sum_{i=1}^N v_i^2 + \frac{\e}{2} \sum_{i \neq j} \varphi \left(\frac {|x_i-x_j|}{\e}\right)\;.
$$
 The corresponding kinetic equation is the Lenard-Balescu equation. 
\end{itemize}
For a recent analysis of these scalings we refer to \cite{NVW1,NVW2}, where allowed classes of interaction potentials are also discussed. Scaling limits and propagation of chaos have a long history starting from the pioneering work of Bogoliubov \cite{Bog}, see e.g.\,\cite{Spbook,PSbonn} for surveys.

The mechanism for which we expect that propagation of chaos actually holds is deeply different in the above cases.
For the mean-field limit,
the force exerted by particle $j$ on particle $i$ is $O(\frac 1N)$ so that the effect of one particle on the other is triviallynegligible. Hence $z_i$ and $z_j$ can be considered as independent random variables in the limit $N\to \infty$, if they were so at time zero. 
 At variance, in the case of the Boltzmann-Grad limit, the forces are strong; but due to the low density regime, the probability of two {\em given} particles  colliding in $(0,t)$ is small, although any given particle is subject to a collision with strictly positive  rate. 
Finally in the weak-coupling limit, the forces are of order  $\frac 1 {\sqrt{\e}} \sim N^{\frac{2}{3}}$, but the potential range is small and the interaction takes place in a time interval of order $\e$.  Overall the variation of momentum of particle $i$ due to the interaction with particle $j$  is of order $\sqrt{\e}$, thus negligible. A similar argument holds in the high-density case.

%
%

Going further in this analysis, one may ask about quantitative estimates on the defect of chaos.
In recent times, the authors have approached this problem in \cite{PPS} by the method of ``correlation errors'' (also called ``$v$-functions''),  which goes back to previous work in the context of kinetic limits \cite{BDP}-\cite{DPTV}, \cite{PS}.
Correlation errors are a tool to measure the tendency of $f^N_j$ to factorize and converge in the mean field limit. 
We will give a precise definition in the next section (Definition \ref{def:correrr} below; see also Def.\,2.1 in \cite{PPS}).
In \cite{PPS}, they have been used to prove a result on the Kac's model, which we recall next.

We consider $N$ particles following a stochastic dynamics: to  each particle, say particle $i$, we associate a velocity $v_i \in \R^3$ \newcommand{\calv}{V} and the vector $$\calv_N=\{v_1, \cdots, v_N\}$$ changes by means of two-body collisions at random times, with random scattering angle.
The probability density $W^N(\calv_N,t)$ evolves according to the master equation (forward Kolmogorov equation) 
\be
\label{K}
\pa_t W^N=\frac 1N T_N W^N
\ee
where
$$
T_N W^N=\sum_{i<j} T_{i,j} W^N ,
$$
\be
\label{eq:defTij}
T_{i,j} W^N (\calv_N)= \int d\omega B(\omega; v_i-v_j)\{W^N(\calv_N^{i,j})-W^N(\calv_N)\}\;,
\ee
and where $\calv_N^{i,j}=\{v_1, \cdots,v_{i-1}, v_i', v_{i+1}, \cdots, v_{j-1},v_j', v_{j+1}, \cdots, v_N\}$ and the pair $v_i' ,v_j'$ gives the outgoing velocities after a collision with scattering (unit) vector $\o$ and incoming velocities $v_i, v_j $. $\tfrac{B(\omega; v_i-v_j)}{|v_i-v_j|}$ is the differential cross-section of the two-body process and we shall assume here, for simplicity, that $B$ is bounded.
The resulting kinetic equation is the homogeneous Boltzmann equation
\bea
\label{Beq}
\pa_t f(v) &=&  \int dv_1 \int d\omega B(\omega; v-v_1) \{ f(v') f(v_1') -f(v) f(v_1) \} \nn \\
&=:& Q(f,f)(v)
\eea
 Such a model has been introduced and studied by Kac in \cite{Kac,Kac2} and has given rise to a wide literature later on (see e.g.\,\cite{Gr71,Sz91,GM97,HaMi,MM,MMW}). For instance in \cite{PPS}, we prove the following.
\begin{Thm}[\cite{PPS}]
\label{mainPPS}
For all $t > 0$ and all $j = 1, . . . , N$, the marginals satisfy 
\be \label{eq:ThmPPS}
\|f^N_j (t) - f(t)^{\otimes j}\|_{L_1} \leq C_2\, e^{C_1t}\, \frac{j^2}{N}
\ee
where $C_1,C_2>0$ are explicitly computable constants.
\end{Thm}

In the present paper, we are going to recover a similar result, albeit in a different setting which we refer to as {\em grand canonical ensemble}.  This is a standard  formalism in statistical mechanics, adapted here to the kinetic problem. With respect to the model introduced above, the novelty of the grand canonical framework is twofold: (i) we allow the total number of particles $N$ to be a (Poisson) random variable with intensity $\mu \to \infty$, and consequently (ii) we do not label the particles from $1$ to $N$ (as done in the definition of marginal), but compute only averages of symmetrized quantities. This leads to the notion of {\em correlation functions}, denoted below $\left\{ f^\m_j \right\}_{j = 1}^{\infty}$, which play the role of the marginals $\left\{ f^N_j \right\}_{j = 1}^{\infty}$ in this context.

The advantage of the grand-canonical formalism is that the correlations due to conservation of energy and mass are zero by assumption, and one can focus more easily on the correlations having a purely dynamical origin. In this spirit
we will reformulate the method of \cite{PPS} leading to a result analogous to \eqref{eq:ThmPPS}, in terms of correlation functions (Theorem \ref{main} below). We shall see that this leads to simplifications of the corresponding proof in the canonical ensemble. Indeed, even if  $f^N_j$ and $f^\m_j$ are asymptotically equivalent, the latter functions encode fewer correlations. This is reflected eventually  in a smaller number of terms in the basic evolution equations for the error (see Remark \ref{rmk:simpeq} below).

The paper is organized as follows. In the next section, we introduce the Kac's model in the grand canonical ensemble, and define the correlation errors. In Section \ref{sec:GC} we discuss the evolution equations, whose derivation is postponed to the appendix.  In Section 3 we establish the main result, Theorem \ref{main}, and its proof is presented in the last section.

\section{Kac model and grand canonical ensemble} \label{sec:GC}

We start by recalling the main ensembles which can be adopted to describe a collisional dynamics, in the spirit of the abovementioned Kac's model.

\begin{itemize}

\item {\em Microcanonical.} The binary collisions preserve energy, and the stochastic dynamics lives in the energy surface
$$
E=\frac 1N \sum_i v_i^2\;.
$$

\item  {\em Canonical.} $E$ is random and only $N$ is fixed, as in the model discussed in the Introduction. 

\item {\em Grand canonical.} $N$ is random. We introduce a sequence of initial density distributions $\{ W_0^n \}_{n\geq 0}$ on the phase space 
$\bigcup_{n\geq 0}   \R^{3n}$. By definition, $\frac 1{n!} W^n_0(\calv_n) $ is the symmetric probability density of the configuration $\calv_n$, being $n$ the total number of particles. The average number of particles is then
$$
\mu := \sum_{n \geq 0} n \int \frac 1{n!} W^n_0(\calv_n) d\calv_n\;.
$$
Here $\mu >0$ is a free parameter (eventually $\mu \to \infty$) and the distributions do depend on $\m$, although this dependence is dropped in our notation  ($W^n_0 = W^{n,\mu}_0$).
At time $t>0$, the state of the system is described by the time-evolved measure  $\{ W^n (t) \}_{n\geq 0}$ solving (cf.\,\eqref{K})
\be
\label{KgrandC}
\pa_t W^n=\frac 1\m T_n W^n\;.
\ee
\end{itemize}

Focusing now on the grand canonical case, let us introduce a natural description of correlations and propagation of chaos. In order to examine correlations due exclusively to the dynamics, we choose a perfectly tensorized and still symmetric initial state, namely
\be \label{eq:initialdatagc}
W_0^n=f_0^{\otimes n} e^{-\mu} \mu^n
\ee
(however most of our discussion could be extended to a larger class of initial measures).

Instead of computing probability densities of particles with labels $\{1,\cdots,j\}$, 
we define rescaled correlation functions by
$$
f^\mu_j (\calv_j,t)= \mu^{-j} \sum_{k \geq 0} \frac 1 {k!} \int dv_{j+1} \cdots dv_{j+k} W^{j+k} ( \calv_{j+k} ,t )\;.
$$
Loosely speaking, this corresponds to computing the amount of $j-$tuples of particles in the configuration $\calv_j$ at time $t$. For a purely factorized state the average total number of $j-$tuples is $\m^j$, which explains the rescaling factor. In particular with the choice \eqref{eq:initialdatagc}, one checks that $\| f^\mu_j \|_{L^1} =1$ (initially and for all times).

Suppose that
\be
\label{pc1}
\lim_{N \to \infty} f^\m_j =  f(t)^{\otimes j}\;,
\ee
say in $L^1-$norm. We may ask about convergence rates. We observe that, in general (e.g.\,for $j \sim \mu$) $f^\m_j(t)$ is not expected to be asymptotically equivalent to  $(f_1^\m(t))^{\otimes j}$, even if it is assumed to be so at time zero. One can then pose the question: how large can  $j=j(\m)$ be in such a way that \eqref{pc1} holds true?

Having at our disposal the estimate
\be
\label{est1}
\| f^\m_j -   f(t)^{\otimes j} \|_{L^1} \leq C^j \left(\frac 1 \m \right)^ \a
\ee
for some $\a \in (0,1)$ and $C>0$, we cannot make better predictions than $j(\m) \sim \log \m$.
On the other hand if for a moderately large $j$, say $j\approx \mu^\gamma$, $\gamma \in (0,1)$, $f^\m_j \approx (f_1^\m(t))^{\otimes j}$, then it is natural to consider
\be
\label{prod}
 (f_1^\m(t)-f(t) )^{\otimes j} =: E_j\;;
\ee
namely the product of the differences rather than the difference of the products. Expanding \eqref{prod} one finds that
\be\label{cumm}
E_j (t) := \sum_{K \subset J} (-1)^{|K|}  f^\m_{J \backslash  K} (t) f(t)^{\otimes K} 
\ee
where $J=\{ 1,2, \cdots, j \}$ is the set of the first $j$ indices, $K$ is any subset of $J$, $J \backslash K$ 
is the relative complement of $K$ in $J$ and $|K| $ is the cardinality of $K $.
$f_A ^\m(t)$ stands for the $|A|-$marginal $f_{|A|} ^\m(t)$ computed in the configuration $ \left(z_i \right)_{i\in A} $. 
Similarly $f(t)^{\otimes K} = f(t)^{\otimes |K|}$ evaluated in $ \left(z_i \right)_{i\in K}$.
The formula \eqref{cumm} has an inverse formula, (see \cite{PPS}) that is
\be
\label{E2}
f^\m_j (t) = \sum_{K \subset J}  E_{J \backslash K}(t)  f(t)^{\otimes K}
\ee
where the notation $E_{J \backslash K}$ is the same one as for $F^\m_{J \backslash K}$,
and where we set $E_0 =E_{\emptyset}=1$.

This motivates the following
\begin{defi}[Correlation error (grand-canonical)]\label{def:correrr}
For any $j \geq 1$, setting $J=\{ 1,2, \cdots, j \}$, we
define the ``correlation error'' of order $j$ by
\be
\label{E1}
E_j (t) := \sum_{K \subset J} (-1)^{|K|}  f^\m_{J \backslash  K} (t) f(t)^{\otimes K}\;,
\ee
where the terms $K = \emptyset$ and $K = J$ have to be interpreted as $f^\m_J = f^\m_j$
and $(-1)^j f^{\otimes J}=(-1)^j f^{\otimes j}$, respectively.
\end{defi}
For instance 
\bea
E_1 (v_1)&=&f^\m_1(v_1)-f(v_1)  \nn \\
E_2(v_1,v_2)&=& f_2^\m(v_1,v_2)- f_1^\m(v_1) f_1(v_2)-f(v_1) f_1^\m(v_1)+f(v_1) f(v_2)\nn\\
E_3 (v_1,v_2,v_3)&=& f_3^\m(v_1,v_2, v_3) 
-f_2^\m(v_1,v_2) f(v_3)\cdots 
+ f^\m_1(v_1) f(v_2) f(v_3)\cdots\nn\\
&&- f(v_1) f(v_2) f(v_3)\;. \nn
\eea


\section{Dynamical equations} \label{sec:GC}

It is well known, and can be easily checked from \eqref{KgrandC},  that the sequence $\{f^\mu_j \}$ satisfies
 the following hierarchy of equations (called BBKGY in analogy with the Hamiltonian case):
\be\label{eqhiera1}
\pa_t f^\mu_j = \frac {T_j}{\mu} f^\mu_j +  C_{j+1}f^\mu_{j+1}
\ee
for $j=1,\cdots,\infty$\;, where
\be
\label{C}
C_{j+1}f^\m_{j+1} (\calv_j)=\sum_{i=1} ^jC_{i,j+1} f^\m_{j+1} (\calv_j)
\ee
and
\be
C_{i,j+1}f^\m_{j+1}(\calv_j) =\int dv_{j+1} \int d\omega B(\omega; v_i-v_j)\{f^\m_{j+1} (\calv_{j+1}^{i,j+1})-f^\m_{j+1}(\calv_{j+1})\}
\label{eq:Ci}
\ee
with the same notations as in \eqref{eq:defTij}.

This hierarchy must be compared with the limiting hierarchy (in the formal limit $\m\to \infty$) 
\be\label{Bhiera}
\pa_t f_j = C_{j+1}f_{j+1}
\ee
which is satisfied by products $f_j(t)=f(t)^{\otimes j}$,  being $f(t)$ solution of the homogeneous Boltzmann equation
\eqref{Beq}.

Introducing correlation errors given by \eqref {E1}, we deduce from \eqref{eqhiera1} the corresponding evolution equations.
We find (see the Appendix)
\bcrr{
\bea\label{eqhieraerror}
\pa_t E_j&=& \frac{T_j}{\mu}E_j +
D_j E_j + 
D_j^{1}
\left(E_{j+1}\right)  
+D_j^{-1}
\left(E_{j-1}\right) +
D_j^{-2}
\left(E_{j-2}\right)
\eea
}
where $ D_j^1=C_{j+1}$ and the operators $D_j , D_j^{-1},  D_j^{-2}$ are defined by:
\bea
&& D_j E_j = \sum_{i\in J} C_{i,j+1}[ f^{\otimes \{i\}} E_{J^i \cup \{j+1\} } + f^{\otimes \{j+1\}} E_{J } ]\nn\\
&& D^{-1}_j E_{j -1} =  \frac 1\mu \sum_{i,s \in J} T_{i,s}  f^{\otimes \{i\}} E_{J^i } \nn\\
&& D^{-2}_j E_{j-2}  = \frac 1{2\mu} \sum_{i,s \in J} T_{i,s}   f^{\otimes \{i,s\}} E_{J^{i,s} }\;, \label{eq:defsDops}
\eea
where $J$ is the set of the first $j$ indices, $J^i=J \backslash \{i\}$, $J^{i,s} = J \backslash \{i,s\}$ and
$$
f^{\otimes \{i\}}=f(v_i), \,\,\,\, f^{\otimes \{i,s\}}=f(v_i) f(v_s)\;.
$$

Here we use the convention the $ D_1^{-2} =0$. Moreover these operators may depend on time, although we will drop out this dependence from the notation.

Note that the first line contains  operators which do not change the particle number, 
$C_{j+1} $ is an operator decreasing
by one the number of particles, while $D^{-1}_j$ and $D^{-2}_j$ are operators increasing  the number of particles by one and two respectively.

\begin{Rmk} \label{rmk:simpeq}
The marginals in the canonical setting satisfy the same hierarchy of equations as in \eqref{eqhiera1} (with $\m$ replaced by $N$), up to a  combinatorial factor $(N-j)/N$ in front of the collision operator $C_{j+1}$. This slight difference produces several spurious terms in the definitions of the operators $D_j , D_j^{-1},  D_j^{-2}$ (see Eq.s (42)-(44) in \cite{PPS}), which are absent in the definition \eqref{eq:defsDops}. Therefore the grand-canonical equations single out the contributions having a purely dynamical interpretation, in terms of variation of clusters (groups of particles) with $j$ mutually correlated particles. Since the interactions are binary, such a cluster can be created in five ways only, corresponding respectively to the five terms in \eqref{eqhieraerror}, and depending on how many independent particles enter the game through the described collision.
\end{Rmk}


%

We conclude this section by establishing our main result. 

\begin{Thm}
\label{main}
Let us assume $f^\mu_j(0)=f_0^{\otimes j} $ being $f_0$ the initial datum for the Boltzmann equation. Then for all $t>0$, there exists $G>1$ such that, for
$\frac j {\sqrt {\mu} }$ sufficiently small,   
\be
\label{maine}
\| f^\mu_j(t) - f^{\otimes j} (t) \|_{L_1} \leq \frac {Gj^2}{\mu}.
\ee
\end{Thm}

\bigskip

Here we are assuming a full factorization at time zero. Looking at the proof it will be clear that this condition can be relaxed. Here we prefer to work in the simplest context to outline the main ideas.

Let us remark that, as it can be shown from  of the proof of Theorem \ref{main} the condition ``$\frac j {\sqrt {\mu} }$ sufficiently small" can  be quantified as ``$\frac j {\sqrt {\mu} }\leq \alpha e^{-\beta t}$ for some contants $\alpha, \beta>0$".

\section{Towards the Boltzmann equation}

Starting from \eqref{eqhiera1},  by the Dyson (Duhamel) expansion we have:

\begin{eqnarray}
\label{expN}
&&f^\mu_j(t)=
\sum_{n= 0}^{\infty}  
\int_0^t dt_1 \int_0^{t_1} dt_2 \dots \int_0^{t_{n-1}} dt_n  
\\&&
U^\mu_j (t-t_1)C_{j+1} \dots
U^\mu_{j+n-1} (t_{n-1}-t_n)C_{ j+n} U^\mu_{j+n} (t_n)f _{0, n+j} , \nonumber
\end{eqnarray}
where
$f_{0,j}=f_0^{\otimes j}$
and
$$
U^\mu_j (t)=e^{\frac {T_j}\mu t}.
$$
On the other hand we also have, for $f_j(t) =f(t)^{\otimes j} $, being $f(t)$ the solution to the Boltzmann equation
\begin{eqnarray}
\label{expN1}
&&f_j(t)=
\sum_{n= 0}^{\infty}  
\int_0^t dt_1 \int_0^{t_1} dt_2 \dots \int_0^{t_{n-1}} dt_n  
C_{j+1} \dots C_{j+n} f _{0, n+j} \\
&& = \sum_{n= 0}^{\infty}  \frac {t^n}{n!} C_{j+1} \dots C_{j+n} f _{0, n+j} \nn.
\end{eqnarray}

Assuming that $B$ is bounded then we have
\be\label{normtc}
\hspace{0,3cm}\norm{T_j} \leq  C_1 \frac {j(j-1)}2 ,
\hspace{0,8cm} \norm{C_{j+1}} \leq  C_1 j .
\ee
From now on we denote by $\| \cdot \| $ the  $L^1$ -norm and by $C_i, i=1,2 \cdots$ different numerical (positive) constants.

Then since $ \|U^\mu_j(t) \|=1$ both the right hand side of \eqref {expN} and
\eqref {expN1}  
are bounded by
$$
\sum_{n\geq 0} C_2^n \frac {t^n}{n!} j (j+1) \cdots (j+n-1) 
\leq 2^j \sum_{n\geq 0} (2C_2 t)^n,
$$
which is converging for $t$ small. 
Therefore the rest of the two series is negligible, uniformly in $\mu$, and this allows us to perform the limit $\mu \to \infty$.
Indeed realizing that, for all $j$, 
$$
\lim_{\mu \to \infty} U^\mu_j(t) = I
$$
strongly, we have also the term by term convergence of \eqref{expN} to \eqref{expN1} in $L^1$. Once having the convergence for $t < \frac 1{2C_2}$, we can iterate the procedure since the new initial conditions, $f^\mu_j (t)$ satisfy
$\| f^\mu_j(t) \| =1$. Then we have convergence for any arbitrary time.

What is the {\it size of chaos} in this case?
By size of chaos we mean the maximal $\a $ for which we have convergence for all j such that $\frac j{\mu^\alpha}\to 0$ as $\mu\to\infty$.  Our main Theorem \ref{main}  shows that  the size of chaos is 
$\frac 12$. On the other hand, as we shall see in a moment, we do not expect 
a size of chaos strictly greater than $\frac12$ so 
we believe
that our result can be considered as optimal as regards this feature.

We observe  that 
$$
\lim_{\mu, j \to \infty} U^\mu_j(t) = I
$$
only if $\frac {j^2}\mu \to 0$. Indeed expressing $U^\mu(t) $ in terms of the series expansion
$$
 U^\mu(t)  f_{0,j} =\sum_n  \frac {t^n}{n!} \frac 1 {\mu^n} \sum_{i_1<k_1}  \cdots \sum_{i_n<k_n} T_{i_1,k_1} \cdots T_{i_n,k_n} f_{0,j}
$$
we realize that the $L^1$ norm of the right hand side can be bounded by
$$
\sum_n  \frac {(Ct)^n}{n!} \left(\frac {j(j-1) }{2 \mu} \right)^n.
$$ 
 Therefore we expect that an estimate of the size of chaos cannot be better than $\frac 12 $ and this can be considered as optimal.
 
 Coming back to the dynamical description given by \eqref {eqhieraerror}, 
we estimate  the size of the $L^1$ norms  of the operators appearing in \eqref {eqhieraerror} by

$$
\| D_j\| =O(j),\,\,  \| D_j^1 \| =O(j), \,\, \| D_j^{-1} \|=O\left(\frac {j^2} \mu\right),  \,\,  \| D_j^{-2} \| =O \left(\frac {j^2} \mu\right).
$$

Note that \eqref {eqhieraerror} is inhomogeneous so that it has nontrivial solutions even for initial data $E_j(0)=0, \, j>0 $ (namely when the initial state is tensorized) which is the case we consider here. 

The control of $E$ passes trough a Dyson expansion around   
 the two parameters semigroup $\bar V_j$ generated by 
$ \frac{T_j}{\mu}E_j +D_j E_j $, namely the solution of
\bea 
&\pa_t \bar V_j (t,s) = \left(\frac{T_j}{\mu}+
D_j \right)\bar V_j (t,s) \nn \\
&\bar V_j(s,s)=I .\nn
\eea
We  easily derive the following bound
\be
\label{V}
\| \bar V_j(t,s)\| \leq e^{C \left(\frac {j^2} \mu +j\right) (t-s)}.
\ee

The action of the other operators $ C_{j+1}, D_j^{-1}, D_j^{-2}$ represents a positive, a negative  and a negative double jump respectively in the space of indices.

We recall that initially the distribution factorizes, namely
$$
E_j(0) =\delta_{j,0},
$$
then
\begin{eqnarray}
\label{expE10}
&&E_j(t)=
\sum_{n= 0}^{\infty}  \sum_ {k_1 \cdots k_n}
\int_0^t dt_1 \int_0^{t_1} dt_2 \,\,\, \dots \int_0^{t_{n-1}} dt_n  
\\&&
\bar V_{s_1}(t, t_1)D^{k_1} \dots
\bar V_{s_n}(t_{n-1}, t_n)D^{k_n} \bar V_{s_{n+1}} ( t_n)E _{j_0} (0) .\nonumber
\end{eqnarray}
where $k_i \in \{-1,-2,1\}$. The sequence $ \{ s_1 \cdots s_n\} $ is defined recursively by
$$
s_{r+1}=s_{r} +k_r, \qquad s_1=j.
$$
and
$$
j_0=j+\sum_{i=1}^n k_i
$$
the initial index. Here $D^k=D^k_{s_k}$ and we drop the indices for notational simplicity.

We have to consider the case $j_0=0$ but in the following it will be convenient to consider general $j_0\geq 0$.

Note also that  $E_0=1 $ and the only possible term involving $E_0$ is
$$
D^{-2}_2 E_0= \frac 1\mu T_2 f^{\otimes 2} 
$$
which is perfectly defined in  $L^1(\R^{2d})$ being $f=f(t)$ the solution of the homogeneous Boltzmann equation.

For simplicity and without loss of generality we choose the parameters in such a way that
$$
\| T_j \| \leq \frac {j^2} 2, \,\,\,\, \| D_j\| \leq \frac  j2 ,\,\,  \| D^1_j \|\leq j , \,\, \| D_j^{-1} \| \leq \frac {j^2} \mu ,   \,\,  \| D_j^{-2} \| \leq \frac {j^2} \mu 
$$ 
%
%
%
%
%
%
and we can prove
\begin{Prop}

 Suppose that,
\be\label{A}
\| E_j(0)\| \leq \left(\frac {j^2}\mu \right)^{j/2} C_0^j,\ \mbox{ for all }j \geq 1,
\ee
for  some $C_0\geq 1$.
Then, there exists $t_0$ sufficiently small and $A>0$  (both explicitily computable)  such that for any $ t \leq t_0$ , 
\be\label{A1}
\| E_j(t) \| \leq \left(\frac {j^2}\mu \right)^{j/2} (A C_0)^j,\ \mbox{ for all } j \geq 1.
\ee
\end{Prop}

\bigskip

We first observe  that we can assume $j\leq  \frac 2{C_0} \sqrt{\mu}$ because otherwise
$$
\| E_j(t) \| \leq 2^j \leq 2^j \left(\frac {jC_0} {2 \sqrt{\mu}}\right)^j = \left(\frac {jC_0} { \sqrt{\mu}}\right)^j.
$$
Next from \eqref {expE10} we have
\begin{eqnarray}
\label{expE11}
&&E_j(t)=
\sum_{n= 0}^{M-1}  \sum_ {k_1 \cdots k_n}
\int_0^t dt_1 \int_0^{t_1} dt_2 \dots \int_0^{t_{n-1}} dt_n \nonumber  
\\&&
\bar V_{s_1}(t, t_1)D^{k_1} \dots
\bar V_{s_n}(t_{n-1}, t_n)D^{k_n} \bar V_{s_{n+1}} ( t_n)E _{j_0} (0) + \nonumber \\&&
\sum_ {k_1 \cdots k_{M} }\int_0^t dt_1 \dots \int_0^{t_{M-1}} dt_{M} 
\bar V_{s_1}(t, t_1)D^{k_1} \dots
\bar V_{s_{M}} ( t_{M-1}-t_{M} ) D^{k_{M}} E _{j_0} (t_{M}) \nonumber
\end{eqnarray}
and denote by $\cal {T}^1_j(t) $ and $\cal {T}^2_j(t) $ respectively the two terms in the right hand side of the above expression.

Moreover we fix $M= \frac 2{C_0} \sqrt{\mu}$.  With this choice
$$
j+s_m  \leq j+M \leq \frac 4{C_0} \sqrt{\mu}
$$
so that all the operators  $D_j^{-1} $  and $D_j^{-2} $ appearing in the definition of $\cal {T}^1_j$ and $\cal {T}^2_j$ satisfy the bound
$$
\| T_j \|  \leq j, \,\,\, \| D_j^{-1} \| \leq \frac {j^2} \mu \leq j  ,  \,\,\,\,\, \| D_j^{-2}\|  \leq \frac {j^2} \mu \leq j
$$
provided that $\mu$ is large enough. As a consequence, by \eqref{V}
\be
\label{V1}
\| \bar V_j(t)\| \leq e^{jt}.
\ee

We first estimate $\cal {T}^2_j(t)$. Observe that, for $t_0 \leq 1$ for which $ e^{jt_0} \leq e^j$
$$
\| \bar V_{s_1}(t, t_1)D^{k_1} \dots
\bar V_{s_{M}} ( t_{M-1}-t_{M} ) D^{k_{M}} E _{j_0} (t_{M}) \| \leq e^j e^{Mt} (j+M)^M 2^{j+M}
$$
as follows by the estimates
$$
e^{s_1 (t-t_1)} e^{s_2 (t_1-t_2)} \cdots e^{s_M ( t_{M-1}-t_{M} )} \leq e^{(j+M)t }\leq e^j e^{Mt}
$$
$$
(j+k_1)(j+k_1+k_2) \cdots (j+k_1+k_2+ \cdot k_M) \leq (j+M)^M \leq (2M)^M
$$
$$
\|  E _{j_0} (t_{M})  \| \leq 2^{j+M}.
$$
Therefore, setting $\t =t e^t$, using the bound  $M^M\leq C_1^M M! $ and choosing $\t$ so small that $\t 8C_1\leq \frac 12$, we have
$$
\| T^2_j(t) \| \leq  4^M \frac {\tau^M}{M!} M^M e^j 2^{j+M} \leq (2e)^j (\tfrac 12)^M.
$$
 Finally since
$$
\sup_M M^j (\tfrac 12)^M \leq  (C_2j )^j
$$
 we conclude that
\be
\label{T2}
\| \cal {T}^2_j(t) \|  \leq  \left( \frac {C_3 j }{\sqrt{\mu}} \right)^j.
\ee

To estimate $\cal {T}^1_j(t)$ we split it into two contributions $\cal {T}^{1,<}_j(t)$ and $\cal {T}^{1,>} _j(t)$ 
relative to the events $j_0 \leq j$ and $j_0 > j$ respectively. Then, proceeding as above
$$
\| \cal {T}^{1,<}_j(t) \| \leq \sum_{n \leq M-1} \sum_ {k_1 \cdots k_n} e^j \frac {\tau^n}{n!} (j+n)^n \left( \frac {j^2} \mu \right)^{\frac {j-j_0} 2} \| E_{j_0} \|.
$$
 The contribution $ (  \frac {j^2} \mu \big)^{\frac {j-j_0} 2} $ arises from the fact that to reach $j_0$ from $j$ we need at least $\frac {j-j_0} 2$ negative jumps, each of them produces a factor bounded by $\frac {j^2} \mu$.
 
 Finally using \eqref{A} and the usual arguments we arrive at
 $$
 \| \cal {T}^{1,<}_j(t) \| \leq \sum_{n \geq 0} (C_4 \t)^n e^j \left( \frac {j^2} \mu \right)^{\frac {j-j_0} 2} \left( \frac {C_0^2 j_0^2} \mu \right)^{\frac {j_0} 2} \leq  \left( \frac {C_5 C_0^2 j^2} \mu \right)^{\frac {j} 2} .
 $$
 
 Finally 
 $$
\| \cal {T}^{1,>}_j(t) \| \leq \sum_{n \leq M-1} \sum_ {k_1 \cdots k_n} e^j \frac {\tau^n}{n!} (j+n)^n 
\left( \frac { j_0^2} \mu \right)^{\frac {j_0} 2} C_0^{j_0}.
$$
Setting now $ \ell= \ell (k_1,k_2 \cdots k_n) =j_0-j$ we have

\bea
\big( \frac {(j+\ell)^2} \mu\big)^{\frac {j+\ell}2}C_0^\ell & = & \left( \frac {j^2}\mu \right)^{j/2} \left( \frac { (j+\ell)^2} \mu\right)^{\frac {\ell}2} \left( \frac { j+\ell } {j}\right)^jC_0^\ell \nn\\
& \leq & \left( \frac {j^2}\mu \right)^{j/2} (\frac 2{C_0}) ^\ell e^\ell C_0^\ell \leq \,\,\, \left( \frac {j^2}\mu \right)^{j/2} (2e)^n.
\eea
where, in the last step we used that (remind that both $j$ and $M$ are bounded by $ \frac {2 \sqrt {\mu}} {C_0} $)
$$
\frac {(j+\ell)}{\sqrt {\mu} } \leq \frac {j+M}{\sqrt{\mu}} \leq \frac 4 {C_0}, \,\,\,\,\, \left( \frac { j+\ell } {j}\right)^j \leq e^\ell,
\,\,\,\,\, \ell \leq n.
$$
 Thus, for $\t$ small enough:
 $$
\| \cal {T}^{1,>}_j(t) \| \leq  e^{2j} C_0^j  \sum_{n \geq 0} C_5^n  \tau^n 
\left( \frac { j^2} \mu \right)^{\frac {j} 2} \leq \left( \frac { C_6 C_0^2  j^2} \mu \right)^{\frac {j} 2}
$$
and the proof is achieved.

As a Corollary we can easily prove

\begin{Prop}\label{A1}
 Let us suppose that $E_j(0)=0$ for $j\geq 1$ ($E_{\emptyset}(0)=1$) 
corresponding to a fully factorized initial state. Then
 for all  $t>0$ and all $j> 1$, 
 one has 
\be
\label{eqmain}
\| E_j(t) \| \leq \left( A e^{B t } \right)^j \left(\frac{j}{\sqrt \mu}\right)^j 
\ee
where $B \geq 0$ and  $A\geq 1 $.
\end{Prop}

Note that from the evolution equation
$$
\pa_t E_1=D_1E_1+ C_{1,2} E_2 
$$
and Proposition \ref{A1} we also obtain the estimate
$$
\| E_1(t) \| \leq t A^2 e^{2Bt} \frac {4} \mu
$$
which will be useful in a moment.

Proof of Theorem \ref{main} 

 We find that, for some time-dependent constants $C_7,C_8,G$, using Proposition \ref{A1} and the estimate for $E_1(t)$ above,
\begin{align}
\label{EEE}
\| f^\mu_j(t)-f^{\otimes j}(t) \|  \leq & j \|E_1(t) \|+ \sum_{k=2}^j  \frac {j!}{k! (j-k)!} \| E_K  (t) \|   \nn \\
&\leq  \frac {jC_8} \mu +\sum_{k=2}^j  \frac {j!}{k! (j-k)!} 
\left(\frac {A e^{Bt} k}{\sqrt {\mu}} \right)^k  \nn \\
&\leq  \frac {jC_8} \mu +\sum _{k=2}^j  \left(\frac {C_7 j }{\sqrt {\mu}} \right)^k = \frac {jC_8} \mu+ \left(\frac {C_7 j }{\sqrt {\mu}} \right)^2\,\,\,\, \frac 1 {1- \frac {C_7 j }{\sqrt {\mu}}}  \nn \\
&= G\left(\frac { j }{\sqrt {\mu}}\right)^2 \nn  
\end{align}
for $\mu$ large such that $\frac {C_7 j }{\sqrt {\mu}} <\tfrac12$.


\LARGE
\section*{Appendix. Derivation of \eqref{eqhieraerror}} \label{sec:appendix}

In this appendix, we deduce the correlation error equations starting from the BBGKY hierarchy, \eqref{eqhiera1}, which we rewrite for the reader's convenience:
$$
\pa_t f^\mu_j = \frac {T_j}{\mu} f^\mu_j +  C_{j+1}f^\mu_{j+1}\qquad j \geq 1\,.
$$
We shall see that the computation is drastically simpler with respect to the analogous derivation in the canonical ensemble reported in the Appendix of \cite{PPS} (where additional terms appear due to the canonical constraint).

We recall the definition of correlation error
\be
\label{eq:defEjapp}
E_j  := \sum_{K \subset J} (-1)^{|K|}  f^\m_{J \backslash  K} f^{\otimes K}
\ee
together with  the following notation. For any set of particle indices $K = \{i_1,\cdots,i_k\}$ of cardinality $k = |K|$, we write
\begin{eqnarray}
&& f^\m_K = f^\m_k\left( z_{i_1},\cdots, z_{i_k}\right)\;,\nn\\
&& E_K = E_k\left( z_{i_1},\cdots, z_{i_k}\right)\;,\nn\\
&& f^{\otimes K} = f(z_{i_1})\cdots f(z_{i_k})\nn
\end{eqnarray}
with the conventions $f^\m_{\emptyset} = E_{\emptyset} =  f^{\otimes \emptyset} = 1$.
We also write $Q(f,f)_i = Q(f,f)(v_i)$ for the Boltzmann operator \eqref{Beq} evaluated in particle $i$.

We start by the simple computation for $j=1$:
\begin{eqnarray}
&& \pa_t E_{\{1\}} = \pa_t \left( f^\m_{\{1\}} - f^{\otimes \{1\}}\right)\nn\\
&& \ \ \ \ \ \ \ \ \ = C_2 f^\m_{\{1,2\}} - Q(f,f)_1\nn\\
&& \ \ \ \ \ \ \ \ \ = C_2 \left( f^{\otimes\{1,2\}} + f^{\otimes\{1\}}E_{\{2\}}
+ E_{\{1\}}  f^{\otimes\{2\}} + E_{\{1,2\}}\right) - Q(f,f)_1\nn
\end{eqnarray}
where in the last step we used the inverse formula (cf.\,\eqref{E2}) for $j=2$. Noticing that, by \eqref{C}-\eqref{eq:Ci},
$C_2(f^{\otimes\{1,2\}}) = Q(f,f)_1$, we are left with
$$
\pa_t E_{\{1\}} =  C_2 \left( f^{\otimes\{1\}}E_{\{2\}}
+ E_{\{1\}}  f^{\otimes\{2\}}\right) + C_2 \left(E_{\{1,2\}}\right) 
$$
and the two terms on the right hand side correspond to $D_1E_{\{1\}}$ and $D^1_1(E_{\{1,2\}})$ respectively.
Therefore \eqref{eqhieraerror} is verified for $j=1$.

We assume now that $\left\{E_{k}\right\}_{k \leq j-1}$ satisfy \eqref{eqhieraerror}, and prove the same thing for $k=j$.
Using definition \eqref{eq:defEjapp} and the evolution equations for $f^\m_J$ and $f$, we find that
\begin{eqnarray}
&& \pa_t E_J = \pa_t \left( f^\m_J - \sum_{\substack{K \subset J \\ K \neq \emptyset}}    f^{\otimes K}
\,E_{J \backslash K}\right)\nn\\
&&\ \ \ \ \ \ \ = \frac {T_j}{\mu} f^\m_J +  \sum_{i \in J}C_{i,j+1}f^\mu_{J \cup \{j+1\}} 
- \sum_{\substack{K \subset J \\ K \neq \emptyset}}  f^{\otimes K}\,  \partial_tE_{J \backslash K}\nn\\
&&\ \ \ \ \ \ \ \ \ \ - \sum_{\substack{K \subset J \\ K \neq \emptyset}}\sum_{i \in K}  f^{\otimes K \backslash\{i\}} Q(f,f)_i E_{J \backslash K}\;.
\label{eq:appzero}
\end{eqnarray}
We compute the three terms in the second line, one by one. First, using the inverse formula
and denoting $J^{i} = J \backslash \{i\}$, $J^{i,s} = J \backslash \{i,s\}$ etc.,
\begin{eqnarray}
&& \frac {T_j}{\mu} f^\m_J = \frac {T_j}{\mu} E_J + \sum_{\substack{i,s \in J \\ i \neq s}} \frac{T_{i,s}}{\mu}
\left\{ f^{\otimes \{i\}} E_{J^{i}} + \frac{1}{2} f^{\otimes\{i,s\}}E_{J^{i,s}}\right\}\nn\\
&&\ \ \ \ \ \ \ \ \ \ + \sum_{\substack{K \subset J \\ K \neq \emptyset}} f^{\otimes K} \frac{T_{j-k}}{\m} E_{J \backslash K}\nn\\
&&\ \ \ \ \ \ \ \ \ \  + \sum_{\substack{K \subset J \\ |K| \neq 0,1}}  \sum_{\substack{i \in K \\ s \in J \backslash K}} \frac{T_{i,s}}{\mu}
f^{\otimes\{ i\}} f^{\otimes K^i} E_{J \backslash K}
+ \sum_{\substack{K \subset J \\ |K| \neq 0,1,2}} \sum_{\substack{i,s \in K \\ i \neq s}} 
\frac{T_{i,s}}{2\m} f^{\otimes\{i,s\}} f^{\otimes K^{i,s}}E_{J \backslash K}
\label{eq:appfirst}
\end{eqnarray}
and, similarly,
\begin{eqnarray}
&& \sum_{i \in J}C_{i,j+1}f^\mu_{J \cup \{j+1\}} = \sum_{i \in J} C_{i,j+1}
\left\{ E_J f^{\otimes\{j+1\}} + f^{\otimes\{i\}} E_{J^i \cup \{j+1\}}
+ E_{J \cup \{j+1\}}\right\}\nn\\
&& \ \ \ +   \sum_{\substack{K \subset J \\ |K| \neq 0}} \sum_{\substack{i \in J \backslash K}} C_{i,j+1}f^{\otimes K} f^{\otimes\{j+1\}}E_{J \backslash K}\nn\\
&&\ \ \ 
+\sum_{\substack{K \subset J \\ K \neq \emptyset}}\, \sum_{i \in J \backslash K }C_{i,j+1} f^{\otimes K}
\left\{ f^{\otimes \{j+1\}}E_{J \backslash K} + f^{\otimes \{i\}} E_{(J \cup \{j+1\} \backslash K)^i} 
+ E_{J \cup \{j+1\} \backslash K} \right\}\;.
\label{eq:appsecond}
\end{eqnarray}
Secondly, using the inductive hypothesis and the explicit definition of the operators $D_{j-k} , D^1_{j-k},D_{j-k}^{-1},  D_{j-k}^{-2}$, we obtain that
\begin{eqnarray}
&& - \sum_{\substack{K \subset J \\ K \neq \emptyset}}    f^{\otimes K}
\,\partial_t E_{J \backslash K} = -  \sum_{\substack{K \subset J \\ K \neq \emptyset}}
f^{\otimes K} \frac{T_{j-k}}{\m}  E_{J \backslash K}
\nn\\
&&\ \ \ \ \ \  - \sum_{\substack{K \subset J \\ K \neq \emptyset}}\, \sum_{\substack{i,s \in J \backslash K \\ i \neq s}}
\frac{T_{i,s}}{\m} f^{\otimes K} \left\{ f^{\otimes \{i\}} E_{(J \backslash K)^i} + 
\frac{1}{2} f^{\otimes\{i,s\}}E_{(J \backslash K)^{i,s}}
\right\}\nn\\
&&\ \ \ \ \ \ - \sum_{\substack{K \subset J \\ K \neq \emptyset}}\, \sum_{i \in J \backslash K }C_{i,j+1} f^{\otimes K}
\left\{ f^{\otimes \{j+1\}}E_{J \backslash K} + f^{\otimes \{i\}} E_{(J \cup \{j+1\} \backslash K)^i} 
+ E_{J \cup \{j+1\} \backslash K} \right\}.
\label{eq:appthird}
\end{eqnarray}
We note that the first line in \eqref{eq:appfirst} and the first line in \eqref{eq:appsecond} reproduce the right hand side of the desired equation \eqref{eqhieraerror}. Therefore we need to check that the remaining terms cancel out. This is straightforward:
\begin{itemize}
\item the second line in \eqref{eq:appfirst} cancels with the first term on the right hand side in \eqref{eq:appthird};
\item last line in \eqref{eq:appfirst} cancels with the second line in \eqref{eq:appthird} (change variables $K \to K\cup\{i\}$
and $K \to K \cup \{i,s\}$);
\item second line in \eqref{eq:appsecond} cancels with last line in \eqref{eq:appzero};
\item last line in \eqref{eq:appsecond} cancels with last line in \eqref{eq:appthird}.
\end{itemize}
This concludes the derivation of \eqref{eqhieraerror}.

\section*{Acknowledgment}
We thank  the "Istituto Nazionale di Alta Matematica"(INdAM)  for partial support.


\begin{thebibliography}{99}


\bibitem{Bog}
N.N. Bogolyubov: \textit{Problems of a dynamical theory in Statistical Physics.} Moscow: State Technical Press (1946) in Russian; English translation in Studies in Statistical Mechanics I, edited by J. de Boer and G. E. Uhlenbeck, part A, Amsterdam: North-Holland (1962).

\bibitem{BDP} C. Boldrighini, A. De Masi, A. Pellegrinotti : \textit {Non equilibrium fluctuations in particle systems modelling Reaction-Diffusion equations}.  Stochastic Processes and Appl. \textbf {42} , 1-30 (1992).
 

\bibitem{CP}  S. Caprino,  M. Pulvirenti:   \textit {A cluster expansion approach to a one-dimensional Boltzmann equation: a validity result} Comm. Math. Phys .
\textbf {166}, 3 (1995), 603-631.


\bibitem{CDPP91}
S. Caprino, A. De Masi, E. Presutti, M. Pulvirenti:
\textit{A derivation of the Broadwell equation.}
Comm. Math. Phys. {\bf 135} (1991) 3, 443--465. 

\bibitem {CCG99}
E. Carlen, M-C. Carvalho, E. Gabetta:
\textit {Central limit theorem for Maxwellian molecules and truncation of the wild expansion.}
CPAM, {\bf 53} , (2000) 3,  370-397. 

\bibitem {CCG05}
E. Carlen, M-C. Carvalho, E. Gabetta:
\textit {On the relation between rates of relaxation and convergence of wild sums for solutions of the Kac equation.}
Jour. Funct. An. {\bf 220}, (2005)  2, 362-387.

\bibitem{CPW}
S. Caprino, M. Pulvirenti and W. Wagner: \textit{A particle systems approximating stationary solutions to the Boltzmann equation}
 SIAM J. Math. Anal. \textbf{4} (1998), 913-934.
 
  \bibitem{DP}A.De Masi, E. Presutti: \textit {Mathematical methods for hydrodynamical limits}. Lecture Notes in Mathematics 1501, Springer-Verlag, (1991).

\bibitem{DOPT} A. De Masi, E. Orlandi, E. Presutti, L. Triolo: \textit {Glauber evolution with Kac potentials. I.Mesoscopic and macroscopic limits, interface dynamics}.
Nonlinearity \textbf {7}, 633-696, (1994).

\bibitem{DOPT1} A.De Masi, E. Orlandi, E. Presutti, L. Triolo: \textit {Glauber evolution with Kac potentials. II. Fluctuations}. Nonlinearity \textbf {9}, 27--51, (1996).


\bibitem{DPTV} A.De Masi, E. Presutti, D. Tsagkarogiannis, M.E. Vares: 
\textit {Truncated correlations in the stirring process with births and deaths}.
Electronic Journal of Probability,\textbf {17}, 1-35, (2012).

\bibitem{GM97}
C.\,Graham, S.\,M\'{e}l\'{e}ard: \textit{Stochastic particle approximations for generalized Boltzmann models
and convergence estimates}, Annals of Probability {\bf 25} (1997), 115-132.

\bibitem{Gr71}
F.A.\,Gr\"{u}nbaum: \textit{Propagation of chaos for the Boltzmann equation}, Arch. Rational
Mech. Anal. {\bf 42} (1971) 323-345.

\bibitem{HaMi}
M. Hauray, S.\,Mischler: \textit{On Ka's chaos and related problems}, Journal of Functional Analysis {\bf 266} (2014) 6055-6157.

\bibitem{Kac} M. Kac:
\emph{Foundations of kinetic theory},
Proceedings of the Third Berkeley Symposium on Mathematical Statistics
and Probability, University of California Press, Berkeley and Los Angeles, 1956.

\bibitem{Kac2} M. Kac:
{\em Probability and related topics in physical sciences}, Interscience, London-New York, 1959.

\bibitem{MM}
S.\,Mischler, C.\,Mouhot:
\textit{Kac's program in kinetic theory},
Invent. Math. \textbf{193} (2013), 1-147.

\bibitem{MMW}
S.\,Mischler, C.\,Mouhot, B.\,Wennberg:
\textit{A new approach to quantitative propagation of chaos for drift, diffusion and jump processes},
Probability Theory and Related Fields {\bf 161} (2015), 1-2, p. 1-59. 

\bibitem{NVW1}
A. Nota, J.J.L. Vel\'{a}zquez, R. Winter:
\textit{Interacting particle systems with long-range interactions: scaling
limits and kinetic equations.} 
Atti Accad. Naz. Lincei Rend. Lincei Mat. Appl. {\bf 32}:2 (2021), 335-377.

\bibitem{NVW2}
A. Nota, J.J.L. Vel\'{a}zquez, R. Winter:
\textit{Interacting particle systems with long-range interactions: Approximation by tagged particles in random fields.} 
 Atti Accad. Naz. Lincei Cl. Sci. Fis. Mat. Natur. {\bf 33}:2 (2022), 439-506.



\bibitem{PP} T. Paul, M. Pulvirenti: \textit{Asymptotic expansion of the mean-field approximation}. Discrete and Continuous Dynamical Systems A {\bf 39} (2019) 1891-1921.

\bibitem{PPS} T. Paul,  M. Pulvirenti, S. Simonella: \textit{On the Size of Chaos in the Mean Field Dynamics.}  Arch. Rational Mech. Anal. {\bf 231} (2019) 285--317.

\bibitem{PS} M. Pulvirenti, S. Simonella:
\textit{The Boltzmann-Grad limit of a hard sphere system:
analysis of the correlation error}.
Invent. Math.  \textbf{207}(3) (2017) 1135-1237.

\bibitem{PSbonn}
 M. Pulvirenti, S. Simonella:
\textit{A Brief Introduction to the Scaling Limits and Effective Equations in Kinetic Theory.} In Albi, Merino-Aceituno, Nota and Zanella (eds), Trails in Kinetic Theory, SEMA SIMAI Springer Series {\bf 25} (2021) 183-208.


\bibitem{Spbook} 
H. Spohn: \textit{Large Scale Dynamics of Interacting Particles.} Texts and Monographs in Physics, Springer–Verlag, Heidelberg (1991).

\bibitem{Sz91}
A.-S. Sznitman: \textit{Topics in propagation of chaos.} In: \'{E}cole d’\'{e}t\'{e} de Probabilit\'{e}s de
Saint-Flour XIX 1989, Lect. Notes in Math. {\bf 1464}, Springer, Berlin (1991) 165-251.






\end{thebibliography}
\end{document}